\documentclass{amsart}
\usepackage{graphicx}
\theoremstyle{definition}

\theoremstyle{remark}

\numberwithin{equation}{section}

\begin{document}

\title[Some Ergodic Properties of Invertible CA]
{Some Ergodic Properties of Invertible Cellular Automata}%

\author{Hasan AKIN}%
\address{Department of Mathematics, Arts and Science Faculty,
Harran University, Sanliurfa, 63200, Turkey}%
\email{akinhasan@harran.edu.tr}%


\begin{abstract}
In this paper we consider invertible one-dimensional linear
cellular automata (CA hereafter) defined on a finite alphabet of
cardinality $p^k$, i.e. the maps
$T_{f[l,r]}:\mathbb{Z}^{\mathbb{Z}}_{p^k}\rightarrow\mathbb{Z}^{\mathbb{Z}}_{p^k}$
which are given by $T_{f[l,r]}(x) = (y_n)_{n=-\infty }^{\infty} $,
$y_{n} = f(x_{n+l}, \ldots, x_{n+r})
=\overset{r}{\underset{i=l}{\sum }}\lambda _{i}x_{n+i}(\text{mod}\
p^k)$, $x=(x_n)_{n=-\infty }^{\infty}\in
\mathbb{Z}^{\mathbb{Z}}_{p^k}$ and
$f:\mathbb{Z}^{r-l+1}_{p^k}\rightarrow \mathbb{Z}_{p^k}$, over the
ring $\mathbb{Z}_{p^k}$ $(k \geq 2$ and $p$ is a prime number),
where $gcd(p,\lambda_r)=1$ and  $p| \lambda_i$ for all $i \neq r $
(or $gcd(p, \lambda_l)=1$ and $p|\lambda_i$ for all $i\neq l$).
Under some assumptions we prove that  any right (left)
permutative, invertible one-dimensional linear CA $T_{f[l,r]}$ and
its inverse are strong mixing. We also prove that any right(left)
permutative, invertible one-dimensional linear CA is Bernoulli
automorphism without making use of the natural extension
previously used in the literature.


\end{abstract}
\maketitle
\section{Introduction}
\par Cellular automata (CA for brevity), first introduced by Ulam and von Neumann,
has been systematically studied by Hedlund from purely
mathematical point of view \cite{H}. Hedlund's paper started
investigation of current problems in symbolic dynamics. CA have
been widely investigated in a great number of disciplines (e.g.
mathematics, physics, computer sciences, and etc.), the study of
such dynamics  from the point of view of the ergodic theory has
received remarkable attention in the last few years (\cite{A1},
\cite{A2} \cite{BKM}, \cite{CDM}, \cite{CFMM}, \cite{HMM},
\cite{KA}, \cite{K}, \cite{PY}). The dynamical behavior of
$D$-dimensional linear CA (linear CA) over ring $\mathbb{Z}_{m}$
has been studied in \cite{CFMM}. In \cite{S1}, Shereshevsky has
studied ergodic properties of CA, he has also defined the $n$th
iteration of a permutative cellular automata and shown that if the
local rule $f$ is right (left) permutative, then its $n$th
iteration also is right (left) permutative. Blanchard \emph{et
al.} \cite{BKM} have answered some open questions about the
topological and ergodic dynamics of 1-dimensional CA. Pivato
\cite{P} has characterized the invariant measures of bipermutative
right-sided, nearest neighbor cellular automaton. Host \emph{et
al}. \cite{HMM} have studied the role of uniform Bernoulli measure
in the dynamics of cellular automata of algebraic origin.

\par Ito \emph{et al}. \cite{ION} have characterized the invertible
linear CA in terms of the coefficients of its local rule. Manzini
and Margara \cite{MM} have obtained some necessary and sufficient
conditions  for a CA over $\mathbb{Z}_m$ to be invertible. They
have given an explicit formula for the computation of the inverse
of a $D$-dimensional linear CA. They have applied finite formal
power series (\emph{fps} for brevity) to obtain the inverse of a
$D$-dimensional linear CA. The technique of \emph{fps} is well
known for the study of these problems (see \cite{ION} for
details). In \cite{A1, A2}, the author has studied the topological
entropy of an additive CA by using the \emph{fps}.

\par It is well known that there are several notions of mixing (i.e.
weak mixing,  strong mixing, mildly mixing, harmonically mixing
etc.) of measure-preserving transformations on probability space
in ergodic theory. For example, recently, Pivato and Yassawi
\cite{PY} developed broad sufficient conditions for convergence.
They introduced the concepts of harmonic mixing for measures and
diffusion for a linear CA. It is important to know how these
notions are related with each other. In the last few decades, a
lot of work is devoted to this subject (see., e.g. \cite{K},
\cite{S1} and \cite{W}). Kleveland \cite{K} has proved that if $r
< 0$ or $l> 0$, then $T_{f[l, r]}$ is strong mixing, and some of
these CA's even are $t$-mixing. In \cite{A}, the author has
studied some ergodic properties of 1-dimensional linear CA acting
on the space of all doubly-infinite sequences taking values in
ring ${\mathbb{Z}}_{m}$.

\par Although the 1-dimensional linear CA theory and the ergodic
properties of this linear CA have developed somewhat
independently, there are strong connections between ergodic theory
and CA theory. For the definitions and some properties of
1-dimensional linear CA, we refer the reader to \cite{H}, \cite{K}
and \cite{WI} (see also \cite{S1} for details).

\par In \cite{S1}, Shereshevsky shown that if $f$ is right (left) permutative
and $0\leq l<r$ (resp. $l<r\leq 0$), then the natural extension of
the dynamical system $(\mathbb{Z}^{\mathbb{Z}}_{m}, \mathcal{B},
\mu, T_{f[l, r]})$ is a Bernoulli automorphism and also he proved
that if $r<0$ or $l>0$ and $T_{f[l, r]}$ is surjective, then the
natural extension of the dynamical system
$(\mathbb{Z}^{\mathbb{Z}}_{m}, \mathcal{B}, \mu, T_{f[l, r]})$ is
a $K$-automorphism. Later, in \cite{S2}, Sherehevsky has also
shown that if $f$ is left (right) permutative and $l\neq 0$
($r\neq 0$), then the natural extension of $T_{f[l, r]}$ is a
$K$-automorphism. In \cite{S1, S2}, Shereshevsky has used the
natural extension so as to convert a noninvertible CA into an
invertible dynamical system. Kleveland \cite{K} proved that any
bipermutative CA is a Bernoulli system with respect to the uniform
measure. In general, the technical definitions of Bernoulli and
Kolmogorov automorphisms in ergodic theory are studied for
invertible transformations (see \cite{CFS, S1, S2, K, W} for
details). Thus, in this paper, we shall restrict our attention to
certain invertible 1-dimensional linear CA over the ring
$\mathbb{Z}_m$ ($m\geq 2$).

\par In this paper we are only interested in invertible 1-dimensional
linear CA and some of their ergodic properties. We consider
invertible 1-dimensional linear CA defined on a finite alphabet of
cardinality $p^k$, where $p$ is prime number and $k\geq 2$ is an
positive integer. Without loss of generality, we focus on $k=2$.

\par One of the interesting parts of this paper is the idea of
using the simple characterization of the invertibility of an
invertible 1-dimensional linear CA to prove some strong ergodic
properties. Therefore, we think that our results will also give a
possibility of proving certain ergodic properties for a complete
formal classification of invertible multi-dimensional CA defined
on alphabets of composite cardinality.

\par Under some assumptions we prove that any right (left) permutative,
invertible one-dimensional linear CA $T_{f[l,r]}$ and its inverse
are strong mixing. We also prove that any right (left)
permutative, invertible 1-dimensional linear CA is a Bernoulli
automorphism without making use of the natural extension
previously used in the literature (\cite{CFS, S1, S2}).
\par The rest of this paper is organized as follows: In Section 2, we
give basic definitions and notations. In Section 3, we study the
invertibility and permutativity of 1-dimensional linear CA. In
Section 4 we investigate the strong mixing property of this
invertible 1-dimensional linear CA and its inverse. In Section 5,
we study the Bernoulli automorphism. In Section 6, we conclude by
pointing some further problems.

\section{Preliminaries} Let ${\mathbb{Z}}_{m}=\{0, 1,\ldots, m-1\}$
$(m\geq 2)$ be the ring of the integers modulo $m$ and
$\mathbb{Z}^{\mathbb{Z}}_{m}$ be the space of all doubly-infinite
sequences $x=(x_n)_{n=-\infty }^{\infty}\in
\mathbb{Z}^{\mathbb{Z}}_{m}$ and $x_n\in {\mathbb{Z}}_{m}$. A CA
can be defined as a continuous homomorphism over
$\mathbb{Z}^{\mathbb{Z}}_m$ relative to the product topology. The
shift $\sigma :\mathbb{Z}^{\mathbb{Z}}_{m}\rightarrow
\mathbb{Z}^{\mathbb{Z}}_{m}$ defined by $(\sigma x)_{i} = x_{i+1}$
is a homeomorphism of compact metric space
$\mathbb{Z}^{\mathbb{Z}}_{m}$.
\par  $T:\mathbb{Z}^{\mathbb{Z}}_{m}\rightarrow
\mathbb{Z}^{\mathbb{Z}}_{m}$ is defined by $(Tx)_i = f(
x_{i+l},\ldots, x_{i+r})$, where $f:
\mathbb{Z}^{r-l+1}_{m}\rightarrow \mathbb{Z}_{m}$ is a given local
rule or map. It is well known that $T$ commutes with $\sigma$.
Martin \emph{et al.} \cite{MOW} have defined a local rule $f$ to
be linear if it can be written as
\begin{equation}\label{eq1}
f(x_{l},\ldots, x_{r}) = \overset{r}{\underset{i=l}{\sum }}\lambda
_{i}x_{i}(\text{mod}\ m),
\end{equation}
where at least one among $\lambda_{l},\ldots,\lambda _{r}$ is
nonzero, mod $m$. We consider 1-dimensional linear CA $T_{f[l,
r]}$ determined by the local rule $f$:
\begin{equation}\label{eq2}
(T_{f[l, r]}x)=(y_n)_{n=-\infty} ^{\infty}, y_n=f(x_{n+l},\ldots,
x_{n+r})=\overset{r}{\underset{i=l}{\sum }}\lambda
_{i}x_{n+i}(\text{mod}\ m),
\end{equation}
 where $\lambda _{l},\ldots, \lambda _{r}\in \mathbb{Z}_{m}$.\\

\par We are going to use the notation $T_{f[l, r]}$ for linear
CA-map defined in (2.2) to emphasize the local rule $f$ and the
numbers $l$ and $r$.  The \emph{fps} associated with $f$ given in
(\ref{eq1}) is defined as $F(X)=\overset{r}{\underset{i=l}{\sum
}}\lambda _{i}X^{-i}$ (see
 \cite{CFMM} and  \cite{MM} for details).

\section{Invertible 1-dimensional linear CA and permutativity}

\par In this section, we study the invertibility of a
1-dimensional linear CA generated by a linear local rule with
respect to modulo $m$ $(m \geq  2)$ and we investigate the
relation between the invertibility of the CA generated by a linear
local rule $f$
and the permutativity of the local rule $f$.\\
\par The notion of permutative CA was first introduced by Hedlund
in \cite{H}. If the linear local rule
$f:\mathbb{Z}_{m}^{r-l+1}\to\mathbb{Z}_{m}$ is given in (2.1),
then it is permutative in the $j$th variable if and only if
$gcd(\lambda_j,m)=1$, where $gcd$ denotes the greatest common
divisor. A local rule $f$ is said to be right (respectively, left)
permutative, if $gcd(\lambda_r,m)=1$ (respectively,
$gcd(\lambda_l,m)=1$). It is said that $f$ is bipermutative if it
is both left and right
permutative.\\

\textbf{Example 3.1.} Consider the local rule
$f:\mathbb{Z}_{3}^3\rightarrow \mathbb{Z}_{3}$ given by $f(x_{-1},
x_{0}, x_1) = (2x_{-1}+2x_{0}+x_{1})$(\text{mod}
3), then $f$ is both left and right permutative, that is, bipermutative.\\

\par From \cite{S1}, it is clear that if the local rule
$f:\mathbb{Z}_{m}^{u}\rightarrow \mathbb{Z}_{m}$ is right (left)
permutative, then so is its $n$th iteration $f^n:
\mathbb{Z}_{m}^{n(u-1)+1}\rightarrow \mathbb{Z}_{m}$ for every
integer $n\geq 1$. Also, in \cite{S1}, Shereshevsky has stated
that the $n$th iteration $T^n_{f[l, r]}$ of CA $T_{f[l, r]}$
generated by the additive local rule $f$ coincides with the CA $T_{f^n[nl,\ nr]}$.\\

\par Ito \emph{et al.} \cite{ION} characterize the invertible linear CA in
terms of the coefficients of the local rule. They have shown that
if
$T_{f[l,r]}:\mathbb{Z}_m^{\mathbb{Z}}\to\mathbb{Z}_m^{\mathbb{Z}}$
is the linear CA given by $f$:
\begin{equation*}
T_{f[l,r]}(x)(n)=\sum_{j=l}^{r}\lambda_jx(n+j)\ (\text{mod}\ m),
\end{equation*}
then $T_{f[l,r]}$ is invertible if and only if for each prime
factor $p$ of $m$ there exists a unique coefficient $\lambda_j$
($l\leq j\leq r$) such that gcd($p,\lambda_j)=1$, that is, $p\mid
\lambda_i$ and gcd($p,\lambda_j)=1$ for all $i\neq j$.\\

In this way, if $m=p^k$ with $p$ prime, then $T_{f[l,r]}$ is
invertible and left permutative (respectively, invertible and
right permutative) if and only if $gcd(p,\lambda_l)=1$ and $p\mid
\lambda_i$ for $i\neq l$ (respectively, $gcd(p,\lambda_r)=1$ and $p\mid \lambda_i$ for $i\neq r$).\\

\par In this paper, we only consider the following results
originally stated for higher dimensions in \cite{MM} for 1-dimensional linear CA.\\

\par It ic clear that for $m=p^k$ we can state
that \emph{fps} $F$ is invertible if and only if there exists a
unique coefficient $\lambda_j$ such that $gcd(\lambda_j, p)=1$.
So, \emph{fps} $F$ can be written as follows:
\begin{equation}\label{eq4}
F(X)=\lambda_{j}X^{-j}+pH(X),
\end{equation}
where $gcd(\lambda_{j}, p) = 1$ and
$H(X)=\overset{r}{\underset{i=l,i\neq j }{\sum}}\lambda
_{i}X^{-i}$.\\
\par We need the following Theorem to concentrate on the problem of
inverting a finite \emph{fps}, associated to 1-dimensional linear CA, over a prime power.\\

\textbf{Theorem 3.2.} (\cite{MM}, Theorem 3.2) Let $F(X)$ denote
an invertible finite \emph{fps} over $\mathbb{Z}_{p^{k}}$, and let
$\lambda_j$ and $H$  be defined as in (3.1). Let $\lambda_j^{-1}$
be such that $\lambda_j^{-1}.\lambda_j\equiv 1$ (\text{mod} $p$).
Then, the inverse of \emph{fps} $F$ is given by
\begin{equation}\label{eq5}
G(X)=\lambda_j^{-1}X^{-j}(1+p\widetilde{H}(X)+
p^2\widetilde{H}(X)^2+\ldots+p^{k-1}\widetilde{H}(X)^{k-1}),
\end{equation}
where $\widetilde{H}(X) =-\lambda_j^{-1}X^{-j}H(X)$.\\

It is clear that if $F(X)$ is given as the equation (3.1), then
the local rule $f$ associated with $F(X)$ can be written as
follows;
$$
f(x_l,\ldots, x_r ) =
\lambda_jx_j+p\overset{r}{\underset{i=l,i\neq j }{\sum }}\lambda
_{i}x_{i}( \text{mod} p^{k} ).
$$

\textbf{Example 3.3.} Let $f(x_1, x_2, x_3)= 2x_1 +2x_2 +x_3
(\text{mod}\ 2^{2})$. Then it is easy to see that $f$ is right
permutative but is not left permutative. The finite \emph{fps} $F$
associated with the local rule $f$ is
\begin{eqnarray*}
F(X) &=& 2X^{-1}+2X^{-2}+X^{-3}\\
&=& X^{-3}[1+(2X^{1}+2X^{2})].
\end{eqnarray*}
Thus, we can find the inverse of $F$ as follows:\\
$$G(X) = X^{3}[1-2(X^{1}+X^{2})].$$
The local rule $g$ related to the finite \emph{fps} $G(X)$ is
$$g(x_{-5},\ x_{-4},\ x_{-3})= 2x_{-5}+ 2x_{-4}+ x_{-3}
(\text{mod}\ 2^{2}).
$$
It is clear that $g$ is not left permutative. Thus, we conclude
that if $f$ is left (respectively, right) permutative and
invertible, then its inverse $g$ is also
left (respectively, right) permutative.\\

\par In Section 4 and 5, we will need some of the identities that will
appear in the proof of Proposition 3.4. Hence, we will include
the proof of Proposition 3.4.\\

 \textbf{Proposition 3.4.} Given 1-dimensional linear CA
\begin{equation}\label{eq6}
(T_{f[l, r]}(x))_n =f(x_{n+l},\ldots,
x_{n+r})=\overset{r}{\underset{i=l}{\sum }}\lambda _{i}x_{n+i}\
(\text{mod}\ p^{k}),
\end{equation}
where $p$ is a prime number. If $gcd(p, \lambda_r)=1$ and
$p|\lambda_i$ for all $i\neq r$ (respectively, $gcd(p,
\lambda_l)=1$ and  $p| \lambda_i$ for all $i \neq l $),  then $f$
is right (respectively, left) permutative and $T_{f[l, r]}$ is
invertible.
\begin{proof}
For brevity without loss of generality we focus on $k=2$. For
$k>2$, similarly the proof can be satisfied.\\

Let us consider the following local rule
\begin{equation}\label{eq7}
f(x_l,\ldots, x_r) =p\overset{r-1}{\underset{i=l}{\sum}}\beta_ix_i
+\lambda_rx_r\ (\text{mod}\ p^{2}),
\end{equation}
where $gcd(p, \lambda_r)=1$ and $\beta_i\in \mathbb{Z}_{p^{2}}$.
From Theorem 3.3, it is clear that the inverse of the local rule
(\ref{eq7}) is as follows:
\begin{equation}\label{eq8}
g(x_{-(2r-l)},\ldots, x_{-r})
=-p\overset{-r}{\underset{i=-(2r-l)}{\sum}}\gamma_ix_i
+\lambda_{r}^{-1}x_{-r}\ (\text{mod}\ p^{2}),
\end{equation}
where $\gamma_i\in \mathbb{Z}_{p^{2}}$. Thus, if 1-dimensional CA
generated by the local rule $f$ is $T_{f[l, r]}$, then the inverse
of $T_{f[l, r]}$ is $T_{f[l, r]}^{-1}
=T_{g[-(2r-l), -r]}$, generated by the local rule $g$.\\

\par Similarly, let us consider the following local rule
\begin{equation}\label{eq9}
f(x_{-l},\ldots, x_{-r})
=\lambda_{-l}x_{-l}+p\overset{-r}{\underset{i=-l+1}{\sum}}a_ix_i \
(\text{mod}\ p^{k}),
\end{equation}
where $gcd(p, \lambda_{-l})=1$, $0 < l <r$ and $p$ is a prime
number, for brevity again without loss of generality we focus on
$k=2$. So we can obtain the finite \emph{fps} $F$ associated with
$f$ as follows:
$$
F(X)=\lambda_{-l}X^l(1+p\lambda_{-l}^{-1}\overset{-r}
{\underset{i=-l+1}{\sum}}a_iX^{-i-l}).
$$
From Theorem 3.3, the inverse of $F$ is obtained as
$$
G(X)=\lambda_{-l}^{-1}X^{-l}(1-p\lambda_{-l}^{-1}
\overset{-r}{\underset{i=-l+1}{\sum}}a_iX^{-i-l}).
$$
Therefore, the inverse of local rule $f$ is obtained as follows:
\begin{equation}\label{eq10}
g(x_l,\ldots, x_{(2l-r)})=
(\lambda_{-l}^{-1}x_{l}-p\overset{-r}{\underset{i=-l+1}
{\sum}}(\lambda_{-l}^{-1})^2 a_ix_{2l+i})\ (\text{mod}\ p^{2}).
\end{equation}
\end{proof}

\section{The ergodic properties of invertible 1-dimensional linear CA}
In this section, we study some ergodic properties of the
invertible 1-dimensional linear CA. To be specific, Coven and Paul
\cite{CP} showed that any surjective CA preserves the uniform
measure. This result is restated by Shereshevsky \cite{S1} and
reproved by Kleveland \cite{K}. Next, Kleveland \cite{K}
established that any one-dimensional left- or right-permutative CA
is mixing.\\

\textbf{Definition 4.1.} Let $T:X\rightarrow X$ be a
measure-preserving transformation on a probability space $(X,
\mathcal{B}, \mu)$; then $T$ is called \emph{strong mixing} if for
any $A, B\in \mathcal{B}$ it satisfies the following equation;
\begin{equation}\label{eq11}
\underset{n\rightarrow \infty }{\lim }\mu(T^{-n}A\cap B)=\mu
(A)\mu (B).
\end{equation}
\par In order to prove the main results of this section and Theorem 5.4,
we consider the $\sigma$-algebra $\mathcal{B}$ of Borelian sets of
$\mathbb{Z}_m^{\mathbb{Z}}$ and the uniform Bernoulli probability
measure $\mu:\mathcal{B}\to [0,1]$, which is defined in the
cylinders $C=\mbox{}_a[j_0,j_1,\cdots,j_s]_{s+a}= \{x\in
\mathbb{Z}_m^{\mathbb{Z}}: x_a=j_0, \cdots, x_{a+s}=j_{s}\}$ as
$\mu(C)=m^{-(s+1)}$.\\

We can easily verify that the Bernoulli measure is defined as
follows:
\begin{eqnarray*}
\mu( _{a}[j_{0}, j_{1},\ldots,
j_{s}]_{s+a})&=&\mu(\{x\in \mathbb{Z}_{m}^{\mathbb{Z}}:x_a=j_{0},\ldots, x_{a+s}=j_{s}\})\\
&=& p_{(j_{0})}.p_{(j_{1})}\ldots p_{(j_{s})},
\end{eqnarray*}
where $(p_{(0)}, p_{(1)},\ldots, p_{(m-1)})$ is a probability
vector. If for all $i, j\in \mathbb{Z}_{m}$ equality $p_{(i)}=
p_{(j)}$ holds then $\mu$ is called uniform Bernoulli measure. It
is a general fact that every surjective CA preserves the uniform
Bernoulli measure (see \cite{H, K, S1} for details).\\

\textbf{Proposition 4.2.} Let $T_{f[l, r]}$ be an invertible
1-dimensional linear CA over $\mathbb{Z}_{p^k}$ and
$gcd(p,\lambda_r)=1$, $p|\lambda_i$ for all $i\neq r$. Then,
$T_{f[l, r]}$ is uniform Bernoulli measure-preserving
transformation.
\begin{proof}
Let $T_{f[l, r]}$ be an invertible 1-dimensional CA and consider a
cylinder set
$$
C = _{a}[j_{0}, j_{1},\ldots, j_{s}]_{s+a}=\{x\in
\mathbb{Z}_{p^{k}}^{\mathbb{Z}}:x_a^{(0)}=j_{0},\ldots,
x_{a+s}^{(0)}=j_{s}\}
$$
Then the first preimage of $C$ under $T_{f[l, r]}$ is as follows:
\begin{eqnarray*}
T_{f[l, r]}^{-1}(C)&=&T_{f[l, r]}^{-1}(\{x\in
\mathbb{Z}_{p^{k}}^{\mathbb{Z}}:x_a^{(0)}=j_{0},\ldots,
x_{a+s}^{(0)}=j_{s}\})\\
&=&\overset{}{\underset{i_{0}, i_{1},\ldots, i_{(r-l)+s}\in
\mathbb{Z}_{p^k}^{((r-l)+s+1)}}{\bigcup}}(_{(a+l)}[i_{0},
i_{1},\ldots, i_{(r-l)+s}]_{a+s+r}),
\end{eqnarray*}
where
$x_a^{(0)}=\overset{r}{\underset{i=l}{\sum}}\lambda_ix_{a+i}^{(1)}$
\ (\text{mod} $p^k$) and
$x_{a+s}^{(0)}=\overset{r}{\underset{i=l}{\sum}}\lambda_ix_{r+a+i}^{(1)}$
(\text{mod} $p^{k}$). It is clear that
$$ _{(a+l)}[i_{0},
i_{1},\ldots, i_{(r-l)+s}]_{a+s+r})=\{x\in
\mathbb{Z}_{p^{k}}^{\mathbb{Z}}:x_{a+l}^{(1)}=i_{0},\ldots,
x_{a+s+r}^{(1)}=i_{r-l+s}\}.
$$
 Then we have
\begin{eqnarray*}
\mu(C)&=&\mu(\{x\in
\mathbb{Z}_{p^{k}}^{\mathbb{Z}}:x_a^{(0)}=j_{0},\ldots,
x_{a+s}^{(0)}=j_{s}\})\\
&=&\mu(T_{f[l, r]}^{-1}(\{x\in
\mathbb{Z}_{p^{k}}^{\mathbb{Z}}:x_a^{(0)}=j_{0},\ldots,
x_{a+s}^{(0)}=j_{s}\})\\
&=&\mu(\overset{}{\underset{i_{0},\ i_{1},\ldots, i_{(r-l)+s}\in
\mathbb{Z}_{p^k}^{((r-l)+s+1)}}{\bigcup}}(_{(a+l)}[i_{0},
i_{1},\ldots,
 i_{(r-l)+s}]_{a+s+r}))\\
&=& (p^{k})^{-(s+1)}.
\end{eqnarray*}
\end{proof}
In \cite{K}, for a CA $T_{f[l, r]}$ Kleveland has proved that if
$r<0$ or $l>0$, then $T_{f[l, r]}$ is strong mixing. He has also
proved that if $f$ is permutative in $x_l$ and $l<0$ or if $f$ is
permutative in $x_r$ and $r>0$, then $T_{f[l, r]}$ is strong
mixing (see \cite{K} for details).\\
\par Under some conditions, now by using a simple
characterization of invertibility we are going to prove that the
invertible linear CA $T_{f[l, r]}$ and its
inverse are strong mixing.\\

\textbf{Theorem 4.3.} Let $T_{f[l, r]}$ be an invertible
1-dimensional linear CA over $\mathbb{Z}_{p^{2}}$ and $gcd(p,
\lambda_r)=1$ and $p| \lambda_i$ for all $i\neq r$, $l>0$. Then
both $T_{f[l, r]}$ and $T_{g[-(2r-l),-r]}$ are strong mixing,
where $g$ is the local rule given in (3.5).
\begin{proof}
Let us firstly consider right permutative local rule $f$ defined
as follows:
\begin{equation*}
f(x_{l}, \ldots, x_{r}) =\sum^{r}_{i=l}\lambda_ix_{i} (\text{mod}\
p^{2}),
\end{equation*}
where $\lambda_i\in  \mathbb{Z}_{p^{2}}$ and $0<l<r$.
\par Let
$A=_a[i_0,\ldots, i_{k}]_{k+a}$ and $B=_b[j_0^{(0)},\ldots,
j_{t}^{(0)}]_{t+b}$ be two cylinder sets. Then we can observe that
$$
A\cap T_{f[l, r]}^{-n}B= \bigcup_{x_{k+1},\ldots, x_{nl-1}}
\bigcup_{j_{nl}^{(n)},\ldots,j_{nr}^{(n)}}
(_a[i_0,\ldots,i_{k},x_{k+1},\ldots,x_{nl-1},j_{nl}^{(n)},\ldots,j_{nr+t}^{(n)}]_{c}),
$$
where $f(j_{i+l}^{(n)}, \ldots, j_{i+r}^{(n)})
=\sum^{r}_{u=l}\lambda_uj_{i+u}^{(n)}$(\text{mod}
$p^{2}$)=$j_{i}^{(n-1)}$ and $c=a+b+k+t+nr$. Thus, for all $nl>k$
we have
\begin{eqnarray*}
&&\mu(A\cap T_{f[l, r]}^{-n}B)=\nonumber\\
&=& \mu(\bigcup_{x_{k+1},\ldots,x_{nl-1}}
\bigcup_{j_{nl}^{(n)},\ldots,j_{nr+t}^{(n)}}
(_a[i_0,\ldots,i_{k},x_{k+1},\ldots,x_{nl-1},j_{nl}^{(n)},
\ldots,j_{nr+t}^{(n)}]_{c}))\nonumber\\
&=&\mu(A)(\sum_{x_{k+1},\ldots,x_{nl-1}} p_{(x_{k+1})}\ldots
p_{(x_{nl-1})})\sum_{j_{nl}^{(n)},\ldots,j_{nr+t}^{(n)}}
p_{(j_{nl}^{(n)})}\ldots p_{(j_{nr+t}^{(n)})}\\
&=&\mu(A)\sum_{j_{nl}^{(n)},\ldots,j_{nr+t}^{(n)}}
p_{(j_{nl}^{(n)})}\ldots p_{(j_{nr+t}^{(n)})}\\
&=&\mu(A)\mu(B).
\end{eqnarray*}
Due to $-r<0$, the Theorem 6.2 in \cite{K} is satisfied. Thus, the
linear CA $T_{g[-(2r-l),-r]}$ is strong mixing.
\end{proof}

It is well known that if the measure-preserving transformation of
probability space is strong mixing, then it is both weak mixing
and ergodic.

\section{Bernoulli automorphism}

As mentioned in introduction, our goal in this section is to show
that certain invertible 1-dimensional linear CA is Bernoulli
automorphism, which is strong ergodic property, without making use
of the natural extension which Shereshevsky \cite{CFS, S1, S2} employs in his proofs.\\

\par Firstly, we recall some definitions from the theory of Bernoulli automorphisms
(see \cite{S1} for details).

\textbf{Definition 5.1.} The partitions $\xi =\{C_i\}$ and
$\eta=\{D_j\}$ of the measure space $(\mathbb{Z}^{\mathbb{Z}}_m,
\mathcal{B}, \mu)$ are called $\varepsilon$-independent
($\varepsilon\geq 0$), if
$$
\overset{}{\underset{i,\ j}{\sum}}|\mu(C_i\cap
D_j)-\mu(C_i)\mu(D_j)|\leq \varepsilon.
$$
The partitions are independent, if they are 0-independent. A
partition $\xi =\{C_i\}$ is called Bernoulli for $T$, that is, an
automorphism of the measure space
$(\mathbb{Z}^{\mathbb{Z}}_m,\mathcal{B}, \mu)$, if all its shifts
are pairwise independent. A partition $\xi =\{C_i\}$ is weakly
Bernoulli for $T_{f[l, r]}$, if for every $\varepsilon > 0$ there
exists an integer $N>0$ such that the partitions
$\overset{0}{\underset{k=-n}{\bigvee}}T^k\xi$ and
$\overset{N+n}{\underset{k=N}{\bigvee}}T^k\xi$ are
$\varepsilon$-independent for all $n\geq 0$.\\

The automorphism $(\mathbb{Z}^{\mathbb{Z}}_m, \mathcal{B}, \mu,
T_{f[l, r]})$ is Bernoulli if and only if it has a
generator $\xi$ which is (weakly) Bernoulli for $T$.\\

\par In this section, our purpose is to prove that the
certain invertible 1-dimensional CA $T_{f[l, r]}$ is Bernoulli
automorphism.\\

\textbf{Lemma 5.2.} If the local rules $f$ and $g$ is defined as
(\ref{eq7}) and (\ref{eq8}), respectively, then partitions
$\overset{n}{\underset{k=0}{\bigvee}} T^{-k}_{g[-(2r-l),
-r]}\xi(-i, i)$ and
$\overset{n}{\underset{k=0}{\bigvee}}T^{-k}_{f[l, r]}\xi(-i,
i)$ are $\varepsilon$-independent.\\
Proof. Follow the definitions.\\

\textbf{Lemma 5.3.} If the local rules $f$ and $g$ is defined as
(\ref{eq7}) and (\ref{eq8}), respectively, then the partition
$\xi(-i,i)$ is weakly Bernoulli for an invertible 1-dimensional
linear CA $T_{f[l, r]}$.
\begin{proof} Let us consider the partition $\xi(-i, i)=\overset{i}{\underset{u=-i}
{\bigvee}}\sigma^{-u}\xi$, where $\sigma$ is the shift
transformation and $\xi$ is the zero-time partition of
$\mathbb{Z}^{\mathbb{Z}}_{p^{2}}$, that is, $\xi=\{_{0}[j]:j\in
\mathbb{Z}_{p^{2}}\}$, then we have
\begin{eqnarray}\label{eq12}
\overset{0}{\underset{k=-n}{\bigvee}}T^{k}_{f[l, r]}\xi(-i,
i)&=&\overset{n}{\underset{k=0}{\bigvee}}T^{-k}_{f[l, r]}\xi(-i,\
i)\nonumber\\
&\preceq&\xi(-i, i)\vee \xi(l-i, r+i)\vee \ldots\vee \xi(nl-i,
nr+i).
\end{eqnarray}
Similarly one gets
\begin{eqnarray}\label{eq13}
\overset{N+n}{\underset{k=N}{\bigvee}}T^{k}_{f[l, r]}\xi(-i, i)&
=&\overset{n}{\underset{k=0}{\bigvee}}(T^{-1})^{-(k+N)}_{f[l, r]}\xi(-i, i)\nonumber\\
&=&\overset{n}{\underset{k=0}{\bigvee}}T_{g[-(2r-l), -r]}^{-(k+N)}\xi(-i, i)\nonumber\\
&\preceq &\xi(-(2r-l)N-i, -rN+i)\vee \times\nonumber\\
& & \xi(-(2r-l)(N+1)-i,
-r(N+1)+i)\times\nonumber\\
& &\vee \ldots\vee \xi(-(2r-l)(N+n)-i,-r(N+n)+i).
\end{eqnarray}
Thus, from (\ref{eq12}) and (\ref{eq13}) it is clear that
$$
\overset{n}{\underset{k=0}{\bigvee}}T^{-k}_{f[l,
r]}\xi(-i,i)\preceq \xi(-i, nr+ i)
$$
and
$$
\overset{N+n}{\underset{k=N}{\bigvee}}T^{k}_{f[l,r]}\xi(-i,i)\preceq
\xi(-(2r-l)(N+n)-i,-rN+i).
$$
For any $n >0$ we have
\begin{equation*}
\overset{p^{2nr+4i+2}}{\underset{a=1}{\sum}}
\overset{p^{2((N+n)(2r-l)-Nr+2i+1)}}{\underset{b=1}{\sum}}|\mu(C_a\cap
C_b)-\mu(C_a)\mu(C_b)|<\varepsilon,
\end{equation*}
where $C_a\in \xi(-i, nr+ i)$ and $C_b\in \xi(-(2r-l)(N+n)-i,
-rN+i)$. Thus, for every $i\geq 0$ the partition $\xi(-i,
i)=\overset{i}{\underset{u=-i}{\bigvee}}\sigma ^{-u}\xi$ is weakly
Bernoulli for the automorphism $T_{f[l, r]}$.
\end{proof}
\textbf{Theorem 5.4.} If the conditions of Lemma 5.2 are
satisfied, then the dynamical system
$(\mathbb{Z}^{\mathbb{Z}}_{p^2}, \mathcal{B}, \mu, T_{f[l, r]})$
is a Bernoulli automorphism.
\begin{proof}  From Lemma 5.3, it is clear that the partition $\xi(-i,
i)=\overset{i}{\underset{u=-i}{\bigvee}}\sigma ^{-u}\xi$ is a
generator for $T_{f[l, r]}$, that is
$\overset{\infty}{\underset{n=-\infty}{\bigvee}}T_{f[l,
r]}^n\xi(-i, i)=\varepsilon$ (see \cite{S1}). Thus, $T_{f[l, r]}$
is a Bernoulli automorphism.
\end{proof}
\textbf{Remark 5.5.} Similarly, one can prove that for $k > 2$ the
dynamical system $(\mathbb{Z}^{\mathbb{Z}}_{p^k},\mathcal{B}, \mu,
T_{f[l, r]})$ defined by the local rule in (\ref{eq7}) is a
Bernoulli automorphism. Also one can prove that if the local rule
$f$ is defined as
$$
f(x_{l},\ldots, x_{r}) =(p_1.p_2 \ldots p_h
\overset{r-1}{\underset{i=l}{\sum}}\lambda_ix_i+\lambda_rx_r)
(\text{mod}\ m),
$$
where $m=p_{1}^{k_{1}}p_{2}^{k_{2}}\ldots p_{h}^{k_{h}}$ and for
all $i\in [l, r]\cap \mathbb{Z}$, $\lambda_i \in \mathbb{Z}_{m}$
and $gcd(\lambda_r,\ m)=1$, then the invertible 1-dimensional CA
$T_{f[l, r]}$
generated by the local rule $f$ is a Bernoulli automorphism.\\
\section{Conclusions}
In this paper, our main results are Theorem 4.3 and Theorem 5.4.
\begin{itemize}
    \item Theorem 4.3 states that any right (left) permutative,
    the invertible linear CA $T_{f[l, r]}$ and its inverse are strong mixing.
    \item Theorem 5.4 states that under the conditions of Lemma
    5.2 the invertible 1-dimensional linear CA $T_{f[l, r]}$ is Bernoulli automorphism.
   \end{itemize}
\par One of the interesting parts of this paper is the idea of
using the simple characterization of the invertibility of an
invertible 1-dimensional linear CA to prove some strong ergodic
properties without making use of the natural extension. This
method provides considerable technical simplifications. Therefore,
we think that our results will also give a possibility of proving
certain ergodic properties for a complete formal classification of
invertible multi-dimensional CA defined on alphabets of composite
cardinality (or the other finite rings). In \cite{AI}, Ak\i n and
\c{S}iap have investigated invertible CA over the Galois rings.
Thus, similar computations and explorations of CA's over different
rings remain to be of interest.

\end{document}